\documentclass{IEEEtran4PSCC} 
\author{
  \IEEEauthorblockN{Sungho Shin and Mihai Anitescu}
  \IEEEauthorblockA{
    Mathematics and Computer Science Division\\
    Argonne National Laboratory,
    Lemont, IL, USA\\
    sshin@anl.gov, anitescu@mcs.anl.gov}
  \and
  \IEEEauthorblockN{François Pacaud}
  \IEEEauthorblockA{Centre Automatique et Systèmes\\
    Mines Paris - PSL,
    Paris, France\\
    francois.pacaud@minesparis.psl.eu}
}
\makeatletter
\let\old@ps@headings\ps@headings
\let\old@ps@IEEEtitlepagestyle\ps@IEEEtitlepagestyle
\def\psccfooter#1{%
    \def\ps@headings{%
        \old@ps@headings%
        \def\@oddfoot{\strut\hfill#1\hfill\strut}%
        \def\@evenfoot{\strut\hfill#1\hfill\strut}%
    }%
    \def\ps@IEEEtitlepagestyle{%
        \old@ps@IEEEtitlepagestyle%
        \def\@oddfoot{\strut\hfill#1\hfill\strut}%
        \def\@evenfoot{\strut\hfill#1\hfill\strut}%
    }%
    \ps@headings%
}
\makeatother

\psccfooter{%
        \parbox{\textwidth}{\hrulefill \\ \small{23rd Power Systems Computation Conference} \hfill \begin{minipage}{0.2\textwidth}\centering \vspace*{4pt} \includegraphics[scale=0.06]{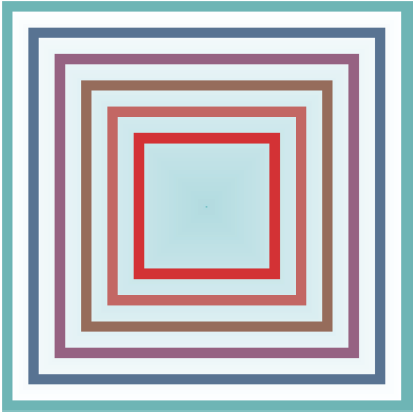}\\\small{PSCC 2024} \end{minipage} \hfill \small{Paris, France --- June 4 -- 7, 2024}}%
}

\usepackage{graphicx}
\usepackage{enumitem}
\usepackage{mathtools}
\usepackage{tikz}
\usepackage{lscape}
\usepackage{multirow}
\usepackage{algorithm}
\usepackage{algorithmic}
\usepackage{cuted}
\usepackage{amsmath,amssymb,amsthm,amsfonts} 
\usepackage[linktocpage=true,colorlinks=true,linkcolor=blue,citecolor=blue,urlcolor=blue]{hyperref}

\newcommand{\st}{\mathop{\text{\normalfont s.t.}}}
\newcommand{\diag}{\mathop{\text{\normalfont diag}}}

\allowdisplaybreaks
\renewcommand{\theenumi}{(\alph{enumi})}

\usepackage{listings}
\lstdefinelanguage{Julia}%
  {morekeywords={abstract,break,case,catch,const,continue,do,else,elseif,%
      end,export,false,for,function,immutable,import,importall,if,in,%
      macro,module,otherwise,quote,return,switch,true,try,type,typealias,%
      using,while},%
   sensitive=true,%
   alsoother={$},%
   morecomment=[l]\#,%
   morecomment=[n]{\#=}{=\#},%
   morestring=[s]{"}{"},%
   morestring=[m]{'}{'},%
}[keywords,comments,strings]%

\usepackage{cite}

\title{Accelerating Optimal Power Flow with GPUs: SIMD Abstraction of Nonlinear Programs and Condensed-Space Interior-Point Methods
}

\begin{document}
\maketitle

\begin{abstract}
  This paper introduces a framework for solving alternating current optimal power flow (ACOPF) problems using graphics processing units (GPUs). While GPUs have demonstrated remarkable performance in various computing domains, their application in ACOPF has been limited due to challenges associated with porting sparse automatic differentiation (AD) and sparse linear solver routines to GPUs. We address these issues with two key strategies. First, we utilize a single-instruction, multiple-data abstraction of nonlinear programs. This approach enables the specification of model equations while preserving their parallelizable structure and, in turn, facilitates the parallel AD implementation. Second, we employ a condensed-space interior-point method (IPM) with an inequality relaxation. This technique involves condensing the Karush--Kuhn--Tucker (KKT) system into a positive definite system. This strategy offers the key advantage of being able to factorize the KKT matrix without numerical pivoting, which has hampered the parallelization of the IPM algorithm. By combining these strategies, we can perform the majority of operations on GPUs while keeping the data residing in the device memory only. Comprehensive numerical benchmark results showcase the advantage of our approach. Remarkably, our implementations—MadNLP.jl and ExaModels.jl—running on NVIDIA GPUs achieve an order of magnitude speedup compared with state-of-the-art tools running on contemporary CPUs.
\end{abstract}


\section{Introduction}

The adoption of graphics processing units (GPUs) in the mathematical programming community has remained
limited due to the challenges associated with parallelizing the optimization algorithms. Notably, nonlinear programming (NLP) remains dependent on
algorithms developed in the 1990s offering limited room for parallelism.
One of the primary challenges arises from the
automatic differentiation (AD) of sparse model equations and the
parallel factorization of indefinite sparse matrices, which are
commonly encountered within constrained numerical optimization tasks
\cite{anitescu2021targeting}. While GPU computation can trivially
accelerate several parts of the optimization process---especially
various internal computations within the optimization solver---the
sluggish data transfer between host and device memory hampers the
ad hoc implementation of GPU accelerations. To
fully leverage the potential offered by modern GPU hardware, a comprehensive computational framework for
optimization on GPUs is imperative. That is, we need an AD/algebraic
modeling framework, sparse linear solvers, and NLP solvers that can
operate entirely on the GPU. Specifically, for the best performance,
both the problem data and the solver's intermediate computational data
must be exclusively resident within the device memory, with the
majority of operations executed on the GPU.

This paper presents our approach to implementing a comprehensive
computational framework for solving large-scale alternating current optimal power flow (ACOPF) problems on
NVIDIA GPUs, along with the associated software implementations:
ExaModels.jl \cite{simdiff}, an algebraic modeling/AD tool, and MadNLP.jl \cite{madnlp}, an
NLP solver. Our approach incorporates two novel strategies: (i) a
single-instruction, multiple-data (SIMD) abstraction of nonlinear
programs, enabling streamlined parallel AD on GPUs, and (ii) a
condensed-space interior-point method (IPM) with an inequality
relaxation strategy, which facilitates the use of highly efficient
{\it refactorization} routines for sparse matrix factorization with
fixed pivot sequences.

While derivative evaluation can be generally cheaper than linear
algebra operations, our numerical results on ACOPF problems show that
AD often constitutes more than half of the total solver time when
using off-the-shelf AD implementations such as
JuMP.jl~\cite{dunning2017jump} or
AMPL~\cite{fourer1990modeling}. Instead, our method leverages a
specialized AD implementation based on the SIMD abstraction of
NLPs. This abstraction allows us to preserve the parallelizable
structure within the model equations, facilitating efficient
derivative evaluations on the GPU. The AC power flow model is
particularly well suited for this abstraction because it involves
repetitive expressions for each component type (e.g., buses, lines,
generators) and the the number of computational patterns does not
increase with the network's size. {\color{black} In other words, the
objective and constraints can be expressed in the form of iterators.}
Numerical results reported in this paper (Table \ref{tbl:results})
demonstrate that our proposed {\color{black} AD strategy} can achieve
over {\color{black} 10} times speedup by running {\color{black} the AD} on the GPU. Compared with general AD
implementations on CPUs (such as AMPL and JuMP.jl), our GPU-based
differentiation method can be approximately 500 times faster.

Linear algebra operations, especially sparse indefinite matrix
factorization, are typically the bottleneck in NLP solution methods.
Parallelizing this operation has been considered challenging,
primarily because of the need for numerical pivoting, which requires irregular
memory accesses and data movement \cite{swirydowicz2022linear}. However, when the matrix can be factorized without
numerical pivoting, a significant part of the operation can be
parallelized, and the numerical factorization can be efficiently
performed on GPUs. We develop a condensed-space IPM strategy that
allows the use of sparse matrix factorization routines without
numerical pivoting. This strategy relaxes equality constraints by
permitting small violations that enable expressing the
Karush--Kuhn--Tucker (KKT) system entirely in the primal space through
the condensation procedure. Although this strategy is not new
\cite{nocedal2006numerical}, it has traditionally been considered less
efficient than the standard full-space method because of increased nonzero
entries in the KKT system. When implemented on GPUs, however, it
offers the key advantage of ensuring positive definiteness in the
condensed KKT system. This, in turn, allows for the utilization
of linear solvers with a fixed numerical pivot sequence (known as
refactorization), an efficient implementation of which is available in the CUDA library.
Although this method is susceptible to numerical stability issues due to
the increased condition number in the KKT system, our results demonstrate
that the solver is robust enough to solve problems with a relative
accuracy of $\epsilon_{\text{mach}}^{1/4}\approx 10^{-4}$.

We present numerical benchmark results to showcase the efficiency of
our method, utilizing our two packages MadNLP.jl and
ExaModels.jl. The  KKT system is solved using the
external cuSOLVER library. To assess the performance of our method, we
compare it with standard CPU approaches using the data available in
pglib-opf \cite{babaeinejadsarookolaee2019power}.  Our benchmark results
demonstrate that our proposed computational framework has significant
potential for accelerating the solution of ACOPF problems, especially
when a moderate accuracy (e.g., $10^{-4}$) is sufficient.
Notably, when running on NVIDIA GPUs, our method achieves a
4x speedup compared with our solver running on CPUs for
the largest instance. Moreover, for the same instance, our approach
surpasses the performance of existing tools (such as Ipopt interfaced
with JuMP.jl) by an order of magnitude. This finding underscores the
importance of harnessing the power of GPUs to tackle
the computational challenges in power systems.

\paragraph*{Contributions}
We present, for the first time, a sparse nonlinear optimization
solution framework that can run entirely on GPUs, with all the
performance-critical data arrays residing exclusively on
GPUs. Additionally, we introduce the concept of SIMD abstraction for
NLP problems, which results in an efficient implementation of
GPU-accelerated parallel AD. Furthermore, we propose the condensed IPM
with an inequality relaxation strategy for the first time, enabling
the treatment of the KKT systems of NLPs without numerical pivoting,
thus allowing the solution of sparse, large-scale NLPs (with a prominent example being ACOPFs) on GPUs.

\paragraph*{Related Work}
Several recent works have explored the use of GPUs for large-scale
nonlinear optimization problems. Cao et al. \cite{cao2016augmented}
proposed an augmented Lagrangian interior-point approach that employs
augmented Lagrangian outer iteration and the treatment of linear
systems using a preconditioned conjugate gradient method. Before the
introduction of the sparse condensed-space IPM with an inequality relaxation
strategy, the authors  investigated the use of reduction
strategies (state variable elimination) to treat KKT systems in a
dense form on the GPU
\cite{pacaud2023parallel,pacaud2022feasible,pacaud2023accelerating}.
In parallel, approaches based on Lagrangian decomposition and
batched with batched TRON solver \cite{lin1999newton}
have been investigated \cite{kim2022accelerated,kim2021leveraging}.
An NLP solver for high-performance computers  with GPU
accelerators, called HiOP, has been under development
\cite{hiop_techrep}, with a scope similar to that of our solver MadNLP.jl.
 Another recent development is the hybrid (direct-iterative) KKT
system solver specifically designed for GPUs \cite{regev2023hykkt}.
The implementation of derivative evaluations with the exploitation of
repeated structures within model equations was, to the best of our
knowledge,  first introduced in Gravity \cite{Gravity}. This was
achieved through the introduction of so-called template constraints,
and multithreaded derivative evaluation has been implemented therein.
We note, however,  that their differentiation approach
is based on symbolic differentiation.
The idea of condensed-space IPM (without inequality relaxation
strategy) is not new, but it has been used primarily in more
specific contexts, where the increased nonzero entries in the KKT system
are less of a concern, such as in the dense form model predictive control
problems \cite{jerez2012sparse,cole2023exploiting}.


\paragraph*{Notation}
We denote the set of real numbers and the set of integers by
$\mathbb{R}$ and $\mathbb{I}$, respectively. We let $[M]:=\{1,2,\cdots,M\}$. We let
$[v_i]_{i\in[M]}:=[v_1;v_2;\cdots,v_M]$.  A vector of ones with an
appropriate size is denoted by $\boldsymbol{1}$. An identity matrix
with an appropriate size is denoted by $I$. For matrices $A$ and $B$,
$A\succ(\succeq) B$ indicates that $A-B$ is positive (semi-)definite
while $A>(\geq) B$ denotes a componentwise inequality. We use the
convention $X:=\diag(x)$ for any symbol $x$.

\section{Preliminaries}\label{sec:prelim}
This section covers three essential background topics: numerical
optimization, AD, and GPU computing.

\subsection{Numerical Optimization}\label{sec:numopt}
We consider NLPs of the following form:
\begin{equation}\label{eqn:cpt}
    \min_{x^\flat \leq x \leq x^\sharp}\;  f(x) \quad \st\;
     g(x) =0 .
\end{equation}
Numerous solution algorithms have been developed in the NLP
literature to solve \eqref{eqn:cpt}. In terms of strategies to deal with inequality
constraints, the NLP solution algorithms can be broadly classified
into active-set methods and interior-point methods \cite{nocedal2006numerical}.
Active-set methods aim to find the set of active constraints associated
with the optimal solution in a combinatorial manner, while IPMs
replace inequality constraints with smooth barrier functions.
IPMs are known to be
more scalable for problems with a large number of constraints and
suitable for parallelization, thanks to the fixed sparsity pattern of
the KKT matrix. Given these advantages, we have chosen IPMs as the
backbone algorithm for developing our optimization methods on GPUs.

In terms of practical computations, three key components play
vital roles: derivative evaluations (often provided by the AD
capabilities of the algebraic modeling languages), linear algebra
operations, and various internal computations within the
solver. Notably, most of the computational efforts are delegated to
the external linear solver and AD library, while the optimization solver
orchestrates the operation of these tools to drive the solution
iterate toward the stationary point of the optimization problem.

Since the successful implementation of the open-source IPM solver
Ipopt, many subsequent implementations of NLP solvers
\cite{chiang2014structured,rodriguez2023scalable,shin2021graph} have
based their implementation on Ipopt
\cite{wachter2006implementation}. We also use Ipopt as our main
reference for the IPM implementation. Below, we outline the overall
computational procedure employed within the NLP solution frameworks.

(1) Given the current primal-dual iterate $(x^{(\ell)},y^{(\ell)},
z^{\flat(\ell)},z^{\sharp(\ell)})$, the AD package evaluates the first-
and second-order derivatives:
\begin{align*}
  \nabla_x f(x^{(\ell)}),\quad
  \nabla_x g(x^{(\ell)}),\quad
  \nabla^2_{xx} \mathcal{L}(x^{(\ell)},y^{(\ell)},z^{\flat(\ell)},z^{\sharp(\ell)}),
\end{align*}
where
\begin{align*}
  \mathcal{L}(x,y,z^{\flat},z^{\sharp}):=f(x) - y^\top
  g(x) - z^\flat (x-x^\flat) - z^\sharp (x^\sharp-x).
\end{align*}

(2) The following sparse indefinite system (known as the KKT system) is
solved using sparse indefinite factorization (typically, via sparse
LBL$^\top$ factorization) and triangular solve routines:
\begin{align}\label{eqn:kkt-indefinite}
  &\begin{bmatrix}
    W^{(\ell)}  + \Sigma^{(\ell)} + \delta^{(\ell)}_w I& A^{(\ell)\top}\\
    A^{(\ell)} & -\delta_c^{(\ell)} I\\
  \end{bmatrix}
  \begin{bmatrix}
    \Delta x\\
    \Delta y\\
  \end{bmatrix}=
  \begin{bmatrix}
    r_x^{(\ell)}\\
    r_y^{(\ell)}\\
  \end{bmatrix},
\end{align}
where
\begin{align*}
  W^{(\ell)}
  &:=\nabla^{2}_{xx}\mathcal{L}(x^{(\ell)},y^{(\ell)},z^{\flat(\ell)},z^{\sharp(\ell)}),
  &&A^{(\ell)}:= \nabla_xg(x^{(\ell)})\\
  \Sigma^{(\ell)}&:= (X^{(\ell)})^{-1}Z^{(\ell)}\\
  r_x^{(\ell)}
  &:=\nabla_x f(x^{(\ell)}) - \mu (X^{(\ell)})^{-1} \boldsymbol{1},
  &&r_y^{(\ell)}:=g(x^{(\ell)}),
\end{align*}
and $\delta^{(\ell)}_w, \delta^{(\ell)}_c>0$ are the regularization parameters
determined based on the inertia correction procedure.

(3) The optimization solver employs a filter line search procedure to
determine the step size~\cite{wachter2006implementation}. The primal-dual iterate is updated by
applying the determined step size and direction. This process is
repeated until the satisfaction of the convergence criteria (typically based on the residual
to the first-order optimality conditions).

\subsection{Automatic Differentiation}
Numerical differentiation of computer programs can be achieved through
three different methods:  finite difference,  symbolic
differentiation, and  AD. The finite difference method suffers
from numerical rounding errors, and its computational complexity grows
unfavorably with respect to the number of function arguments, making
it less preferable unless no other alternatives are
available. Symbolic differentiation uses computer algebra systems to
obtain symbolic expressions of first or higher-order
derivatives. While this method can differentiate functions up to high
numerical precision, it suffers from "expression swelling"  and
struggles to compute the derivatives of long nested expressions in a
computationally efficient way.

In contrast, AD differentiates computer
programs directly by inspecting the computation graph and applying chain rules,
to evaluate derivatives efficiently and accurately.
This approach has become the dominant paradigm for derivative computation within the
scientific computing domain, including NLP and
machine learning. For large-scale optimization problems, such as ACOPFs, AD tools are often implemented as part of domain-specific
modeling languages. Examples of such modeling languages include JuMP,
CasADi, and AMPL (optimization) and TensorFlow, Torch, and Flux
(machine learning).

Derivatives can be propagated through the
recursive application of chain rules in two ways, forward mode and 
reverse mode, which operate in opposite directions (respectively, from leaves to
root and from root to leaves). Reverse-mode automatic
differentiation, also known as the adjoint method, has proven to be
particularly effective for dealing with function expressions in
large-scale optimization problems.

The Julia language, our language of choice, offers convenient and
efficient ways to implement automatic differentiation. Through the use
of the multiple dispatch paradigm~\cite{bezanson2017julia} \textit{any
Julia function}---including commonly used operations such as addition,
multiplication, trigonometric and exponential functions---can be easily overloaded. Multiple dispatch allows functions
to be dynamically dispatched based on the runtime type, a crucial
feature for implementing differentiable programming. Several AD
implementations have been developed in the Julia language, such as
ReverseDiff.jl, ForwardDiff.jl, Zygote.jl, and JuMP.jl. While these
tools are general and useful for various applications, however, they are not
optimized for evaluating derivatives of ACOPF problems, since they are
not designed to exploit the parallelizable structures in the model,
while preserving the desired sparsity.

\subsection{GPU Computing}\label{sec:gpu}
With the increasing prevalence of GPUs in various scientific computing
domains, there has been growing interest in leveraging these emerging
architectures to efficiently solve large-scale NLPs, such as ACOPF problems.
However, adapting an NLP solution algorithm,
such as IPM, to GPUs presents challenges due to the fundamental
differences between GPU and CPU programming paradigms. While CPUs
execute a sequence of instructions on a single input (single
instruction, single data, or SISD, in Flynn's taxonomy), GPUs run the
same instruction simultaneously on hundreds of threads using the SIMD
paradigm. SIMD parallelism works well
for algorithms that can be decomposed into simple instructions running
entirely in parallel, but not all algorithms fit this
paradigm. For example, branching in the control flow can hinder
lockstep execution across multiple threads and, in turn, prevent
efficient implementations on GPUs. On the contrary, when the
algorithm's structure allows for efficient parallelization, the SIMD
parallelism in GPUs can offer orders of magnitude speedup.

We highlight that the following, arguably common, computational
patterns are particularly effective when implemented on GPUs:
\begin{align}
  y&\leftarrow \left[\phi(x; q_j)\right]_{j\in [J]}\tag{Pattern 1}\label{eqn:pattern-1}\\
  o&\leftarrow  \mathop{\rm{Op}}_{i\in [I]} \psi(x;p_i) \tag{Pattern 2}\label{eqn:pattern-2}\\
  x&\leftarrow  \chi_{s_1}\circ\cdots\circ \chi_{s_K}(x) \tag{Pattern 3}\label{eqn:pattern-3}
\end{align}
Here, $\psi:\mathbb{R}^{n_x}\times \mathbb{R}^{n_{p}}\rightarrow
\mathbb{R}$ , $\phi:\mathbb{R}^{n_x}\times
\mathbb{R}^{n_{q}}\rightarrow \mathbb{R}$, and
$\chi:\mathbb{R}^{n_x}\times \mathbb{R}^{n_{s}}\rightarrow
\mathbb{R}^{n_x}$ are simple instructions that require only a small number
of operations; $\mathop{\rm{Op}}$ is a monoid operator on
$\mathbb{R}\cup\{+\infty,-\infty\}$, such as addition, multiplication,
maximum, and minimum. In \ref{eqn:pattern-3}, we denote $\chi_{s_k}(x):=\chi(x,s_k)$ and
assume that $\circ$ is commutative for $\{\chi_{s_k}(\cdot)\}_{\forall
s_k}$. \ref{eqn:pattern-1} is typically most effective on GPUs, where
each thread employed can operate independently without needing to
simultaneously manipulate the same device memory location. While
\ref{eqn:pattern-2} and \ref{eqn:pattern-3} are less effective,
they still can be significantly
faster than operations on CPUs, since a substantial part of the operation
can still be parallelized by the use of buffers. In simple cases, the
implementation of these operations can be performed with the standard
{\tt map} and {\tt mapreduce} programming models. In more
complex cases, however, especially for \ref{eqn:pattern-3}, the implementation
may require preinspection of memory write-access patterns and the use
of custom kernels.

Many of the operations required in AD of sparse physical models, as
well as the application of optimization algorithms, are based on the
computational patterns mentioned above. For example, an ACOPF model
can be implemented with 15 different computational patterns (see
Section~\ref{sec:simd:ad}), all of which fall within the
aforementioned categories. Furthermore, the computation within
optimization solvers, such as forming the left-hand side for the KKT
systems, computing the $\|\cdot\|_\infty$ norm of the constraint
violation, and assembling the condensed KKT system, can be carried out
by using these computational patterns as well. The only exception is
the factorization of the sparse KKT matrix, which requires more
sophisticated implementations.

Implementing kernel functions for the above computational patterns in
the Julia language is straightforward since Julia provides excellent
high-level interfaces for array and kernel programming for GPU
arrays. The code can even be device-agnostic, thanks to the portable
programming capabilities brought by KernelAbstractions.jl. All
of the AD and optimization capabilities in our
tools MadNLP.jl and ExaModels.jl are implemented in  Julia,
by leveraging its kernel and array programming capabilities.


\section{SIMD Abstraction of NLPs}\label{sec:simd}

This section describes our implementation of SIMD abstraction and
sparse AD of the model equations. The abstraction and AD are
implemented as part of our algebraic modeling language ExaModels.jl.

\subsection{Abstraction}
The SIMD abstraction under consideration is as follows:
\begin{align}\label{eqn:prob}
  \min_{x^\flat\leq x \leq x^\sharp}
  & \sum_{l\in[L]}\sum_{i\in [I_l]} f^{(l)}(x; p^{(l)}_i)\\\nonumber
  \st\; &\forall m\in[M]:\\\nonumber
  &\begin{aligned}[t]
    \left[g^{(m)}(x; q_j)\right]_{j\in [J_m]} +\sum_{n\in [N_m]}\sum_{k\in [K_n]}h^{(n)}(x; s^{(n)}_{k}) =0,
  \end{aligned}
  \end{align}
where $f^{(\ell)}(\cdot,\cdot)$, $g^{(m)}(\cdot,\cdot)$, and
$h^{(n)}(\cdot,\cdot)$ are twice differentiable functions with respect
to the first argument, whereas $\{\{p^{(k)}_i\}_{i\in [N_k]}\}_{k\in[K]}$,
$\{\{q^{(k)}_{i}\}_{i\in [M_l]}\}_{m\in[M]}$, and
$\{\{\{s^{(n)}_{k}\}_{k\in[K_n]}\}_{n\in[N_m]}\}_{m\in[M]}$ are
problem data, which can either be discrete or continuous.  We assume
that our functions $f^{(l)}(\cdot,\cdot)$, $g^{(m)}(\cdot,\cdot)$, and
$h^{(n)}(\cdot,\cdot)$ can be expressed with computational
graphs of moderate length. One can observe that the problem in
\eqref{eqn:prob} is expressed by the computational patterns in Section
\ref{sec:gpu}. In particular, the objective function falls within \ref{eqn:pattern-2},
the first term in the constraint falls within \ref{eqn:pattern-1}, and the second term in the constraint falls within \ref{eqn:pattern-3}.
Accordingly, the evaluation and differentiation of the
model equations in \eqref{eqn:prob} are amenable to SIMD parallelism.

To implement the SIMD abstraction in the modeling environment, the
algebraic modeling interface in ExaModels.jl requires the users to
specify the model equations in an {\tt Generator} data type in the
Julia language. This composite data type consists of an instruction (a Julia
function) and data (a host or device array) over which the instruction
is executed. This naturally facilitates maintaining the NLP model
information in the form of SIMD abstraction in \eqref{eqn:prob} and
facilitates the model evaluation and differentiation on GPU
accelerators.

\subsection{Parallel AD}
\label{sec:simd:ad}
Many physics-based models, such as ACOPF, have a highly repetitive
structure. One of the manifestations of it is that the mathematical
statement of the model is concise, even if the practical model may contain
millions of variables and constraints. This is possible due to the use of
repetition over a certain index and data sets. For example,
it suffices to use 15 computational patterns to fully specify the
AC OPF model. These patterns arise from (1) generation cost, (2) reference
bus voltage angle constraint, (3--6) active and reactive power flow (from and to),
(7) voltage angle difference constraint, (8--9) apparent
power flow limits (from and to), (10--11) power balance equations,
(12--13) generators' contributions to the power balance equations, and
(14--15) in/out flow contributions to the power balance
equations. However, such repetitive structure is not well exploited in
the standard NLP modeling paradigms. In fact, without the SIMD
abstraction it is difficult for the AD package to detect the
parallelizable structure within the model, because it will require the full
inspection of the computational graph over all expressions.  By
preserving the repetitive structures in the model, the repetitive
structure can be directly available in the AD implementation.

Using the multiple dispatch feature of Julia, ExaModels.jl generates
highly efficient derivative computation code, specifically compiled
for each computational pattern in the model. These derivative evaluation codes can be run over the data in various GPU array formats
and implemented via array and kernel programming in the Julia Language. In
turn, ExaModels.jl has the capability to efficiently evaluate first- and
second-order derivatives using GPU accelerators.

\subsection{Sparsity Analysis}
The sparsity analysis is needed to determine the sparsity pattern
of the evaluated derivatives. In the case of large-scale sparse
problems, the initial sparsity analysis of nonlinear expressions can
be expensive, since the sparsity should be analyzed for potentially
millions of objective and constraint terms. Often, however, 
these analyses are applied to the same computational patterns, and
the time and memory spent on sparsity analysis can be significantly
reduced if the repetitive structures are exploited.

ExaModels.jl exploits the SIMD abstraction of the model equations to
save the computational cost spent for sparsity analysis.  This is
accomplished by applying sparsity analysis for the instruction for
each computational pattern and expanding the obtained sparsity pattern
over the data over which the instruction is executed. Specifically,
the sparsity analysis code exploits Julia's multiple dispatch feature
to obtain a parameterized sparsity pattern for each instruction, and the
obtained parameterized sparsity pattern is materialized once the data
array is given. This process saves significant computational costs
for the sparsity analysis.

\section{Condensed-Space IPMs with an Inequality Relaxation Strategy}\label{sec:ipm}
We present the condensed-space IPM within the
context of the general NLP formulation in \eqref{eqn:cpt}. Our method
has two key differences from standard IPM
implementations: (i) the use of inequality relaxation and (ii) the
condensed treatment of the KKT system.

\subsection{Inequality Relaxation}

At the beginning of the algorithm, we apply inequality relaxation to replace the equality constraints in \eqref{eqn:cpt} with inequalities by introducing slack variables $s\in\mathbb{R}^{m}$:
\begin{align}\label{eqn:relax}
  g(x)- s = 0,\quad s^{\flat}\leq s\leq  s^\sharp,
\end{align}
where $s^\flat,s^\sharp\in\mathbb{R}^{m}$ are lower and upper bounds chosen to be close to zero.
This relaxed problem can be stated as follows:
\begin{equation}\label{eqn:relaxed}
    \min_{\big[\substack{x^\flat\\s^\flat}\big]\leq \big[\substack{x\\s}\big] \leq \big[\substack{x^\sharp\\s^\sharp}\big]}\;
      f(x)\quad \st\;  g(x) - s = 0.
\end{equation}
In our implementation we set $s^{\flat},s^\sharp$ as
$-\epsilon_{\text{tol}}\boldsymbol{1}$ and $+\epsilon_{\text{tol}}\boldsymbol{1}$, respectively,
where $\epsilon_{\text{tol}}>0$ is a user-specified relative tolerance
of the IPM. This type of relaxation is
commonly used in practical IPM implementations; for example, in Ipopt,
the solver relaxes the bounds and inequality constraints by
$O(\epsilon_{\text{tol}})$ to prevent an empty interior of the
feasible set (see \cite[Section 3.5]{wachter2006implementation}).
For condensed-space IPM, we cannot maintain the same level of precision
because of the increased condition number of the KKT system.
We have found that setting $\epsilon_{\text{tol}}$ to be
$\epsilon_{\text{mach}}^{1/4}\approx 10^{-4}$ ensures numerical
stability while achieving satisfactory convergence behavior. Thus, our
solver sets the $\epsilon_\text{tol}$ to $10^{-4}$ by default when using condensed
IPM.

\subsection{Barrier Subproblem}

The IPM replaces the equality- and inequality-constrained problem in \eqref{eqn:relax} with an equality-constrained  barrier subproblem:
\begin{subequations}\label{eqn:barrier}
  \begin{align}
    \min_{x, s}\;
    &
      \begin{aligned}[t]
        &f(x) - \mu\boldsymbol{1}^\top \log (x- x^\flat) -\mu\boldsymbol{1}^\top \log (x^\sharp -x)\\
        &- \mu\boldsymbol{1}^\top \log (s- s^\flat) -\mu\boldsymbol{1}^\top \log (s^\sharp -s)
      \end{aligned}
          \label{eqn:barrier-obj}\\
    \st\;
    &g(x) - s = 0. \label{eqn:barrier-con}
  \end{align}
\end{subequations}
Here, $\mu>0$ is the barrier parameter. The smooth log-barrier
function is employed to avoid handling inequalities in a combinatorial
fashion (as in active set methods). A superlinear local convergence to
the first-order stationary point can be achieved by repeatedly
applying Newton's step to the KKT conditions of \eqref{eqn:barrier}
with $\mu\searrow 0 $.

\subsection{Newton's Step Computation}
The Newton step direction is computed by solving a so-called KKT system; to explain this, we consider the first-order optimality conditions (KKT conditions) for the barrier subproblem in \eqref{eqn:barrier}:
\begin{align}\label{eqn:first}
  \begin{aligned}[t]
    \nabla_{x} f(x) - \nabla_{x}g(x)^\top y  - z_x^\flat  + z_x^\sharp = 0\;&\\
    \begin{aligned}
      - z_s^\flat  + z_s^\sharp + y &= 0,\\
      Z^\flat_x (x-x^\flat) - \mu\boldsymbol{1} &= 0,\\
      Z^\flat_s (s-s^\flat) - \mu\boldsymbol{1}&= 0,
    \end{aligned}
    \quad
    \begin{aligned}
      g(x) - s  &= 0\\
      Z^\sharp_x (x^\sharp-x) - \mu\boldsymbol{1}&= 0\\
      Z^\sharp_s (s^\sharp-s) - \mu\boldsymbol{1}&= 0,
    \end{aligned}&
  \end{aligned}
\end{align}
where $y\in\mathbb{R}^{m}$, $z_x^\flat,z_x^\sharp\in\mathbb{R}^{n}$,
and $z_s^\flat,z_s^\sharp\in\mathbb{R}^{m}$ are Lagrange multipliers
associated with the equality and bound constraints in
\eqref{eqn:relaxed}. The Newton step for solving the nonlinear
equations in \eqref{eqn:first} can be computed by solving the system in \eqref{fig:very-long-eqn}.
\begin{figure*}[t]
  \vspace{-.2in}
  \begin{align}
    \underbrace{
    \begin{bmatrix}
      W^{(\ell)}  + \delta^{(\ell)}_w I & & A^{(\ell)\top}& -I & I &  \\
      & \delta^{(\ell)}_w I & -I&&&-I & I\\
      A^{(\ell)}& -I & -\delta^{(\ell)}_c I\\
      Z_x^{(\ell)\flat}&&&X^{(\ell)}-X^\flat\\
      -Z_x^{(\ell)\sharp}&&&&X^\sharp-X^{(\ell)}\\
      &Z_s^{(\ell)\flat}&&&&S^{(\ell)}-S^\flat\\
      &-Z_s^{(\ell)\sharp}&&&&&S^\sharp-S^{(\ell)}\\
    \end{bmatrix}
    }_{M_\text{full}}
    \begin{bmatrix}
      \Delta x \\
      \Delta s \\
      \Delta y \\
      \Delta z_x^\flat \\
      \Delta z_x^\sharp \\
      \Delta z_s^\flat \\
      \Delta z_s^\sharp \\
    \end{bmatrix} =
    \begin{bmatrix}
      p^{(\ell)}_{x }\\
      p^{(\ell)}_{s }\\
      p^{(\ell)}_{y }\\
      p^{(\ell)}_{z_x^\flat }\\
      p^{(\ell)}_{z_x^\sharp }\\
      p^{(\ell)}_{z_s^\flat }\\
      p^{(\ell)}_{z_s^\sharp }\\
    \end{bmatrix}
    \label{fig:very-long-eqn} 
  \end{align}
  \vspace{-.3in}
\end{figure*}
Here, we recall the definitions of $W^{(\ell)}$ and $A^{(\ell)}$ from Section \ref{sec:numopt}, and
$p^{(\ell)}_x,\cdots p^{(\ell)}_{z_s^\sharp}$ are defined by the left-hand sides of the equations in
\eqref{eqn:first}. In the sequel, we will drop the superscript
$(\cdot)^{(\ell)}$ for concise notation.

Now, we observe that a significant portion of the system in \eqref{fig:very-long-eqn} can be eliminated by exploiting the block
structure, leading to an equivalent system stated in a smaller
space. In particular, the lower-right $4\times 4$ block is always
invertible since the IPM procedure ensures that the iterates stay in
the strict interior of the feasible set. This allows for eliminating
the lower-right 4x4 block, resulting in
\begin{align}\label{eqn:very-long-reduced}
  &
    \underbrace{
    \begin{bmatrix}
      W  + \Sigma_x + \delta_w I & & A^{\top} \\
      & \Sigma_s + \delta_w I & -I\\
      A& -I & -\delta_c I\\
    \end{bmatrix}
  }_{M_{\text{aug}}}
  \begin{bmatrix}
    \Delta x \\
    \Delta s \\
    \Delta y \\
  \end{bmatrix}=
  \begin{bmatrix}
    q_x \\
    q_s\\
    q_y\\
  \end{bmatrix},
\end{align}
where
\begin{align*}
  \Sigma_x&:= Z^\flat_x (X-X^\flat)^{-1}+ Z^\sharp_x (X^\sharp-X)^{-1}\\
  \Sigma_s&:= Z^\flat_s (S-S^\flat)^{-1}+ Z^\sharp_s (S^\sharp-S)^{-1}\\
  q_x&:=p_x + (X-X^\flat)^{-1} p_{z^\flat_x}-  (X^\sharp-X)^{-1} p_{z^\sharp_x}\\
  q_s&:=p_s + (S-S^\flat)^{-1} p_{z^\flat_s}-  (S^\sharp-S)^{-1} p_{z^\sharp_s}\\
  q_y&:=p_y.
\end{align*}
The bound dual steps can be recovered as follows:
\begin{align}\label{eqn:recover-1}
  \begin{aligned}[t]
    \Delta z^\flat_x &= \left(X-X^\flat\right)^{-1} \left(-Z^\flat_x \Delta x  + p_{z^\flat_x}\right)\\
    \Delta z^\sharp_x &= \left(X^\sharp-X\right)^{-1} \left(Z^\sharp_x \Delta x  + p_{z^\sharp_x}\right)\\
    \Delta z^\flat_s &= \left(S-S^\flat\right)^{-1} \left(-Z^\flat_s \Delta s  + p_{z^\flat_s}\right)\\
    \Delta z^\sharp_s &= \left(S^\sharp-S\right)^{-1} \left(Z^\sharp_s \Delta s  + p_{z^\sharp_s}\right).
  \end{aligned}
\end{align}
Note that the matrices involved in the inversions in \eqref{eqn:recover-1}
are always diagonal, so their computation is cheap.
Also, note that the augmented system in \eqref{eqn:very-long-reduced}
corresponds to the KKT system in \eqref{eqn:kkt-indefinite}. However,
in the original version of the algorithm \cite{pacaud2023accelerating}, we did not introduce the
slack variables, so it did not have the additional structure imposed
by the slack variables.

The key advantage of the inequality relaxation strategy is that it
imposes additional structure on the augmented KKT system, allowing us to
further reduce the dimension of the problem. In particular, the
lower-right 2x2 block in \eqref{eqn:very-long-reduced} can be
eliminated, which is a procedure called {\it condensation}; here, the invertibility of the lower-right block can be verified from the fact that $\delta_w,\delta_c\geq 0$ and $\Sigma_s\succ 0$. Through this,
we obtain the following system, written in the  primal
space only:
\begin{align}\label{eqn:very-long-condensed}
  (\underbrace{W + \delta_wI + \Sigma_x + A^{\top} D A}_{M_\text{cond}} ) \Delta x = q_x + A^\top (C q_s +  Dq_y ),
\end{align}
where
\begin{align*}
  C := \left(\delta_c \Sigma_s + (1+\delta^{}_c\delta^{}_w) I\right)^{-1}, \;
  D := \left(\Sigma_s + \delta^{}_w I\right)C,
\end{align*}
and the dual and slack step directions can be recovered by
\begin{align}
  \Delta s &:= C \left(\delta_c q_s - (q_y + A\Delta x)\right)\nonumber\\
  \Delta y &:= (\Sigma_s + \delta_w I) \Delta s -q_s.\label{eqn:recover-2}
\end{align}
Again, the matrices involved in the inversions above are always
diagonal, so their computation is cheap.

Therefore, the only sparse matrix that needs to be factorized is the
matrix in the left-hand side of \eqref{eqn:very-long-condensed}, with
dimension $n \times n$.  Although we call
\eqref{eqn:very-long-condensed} a condensed KKT system, 
$M_\text{cond}$ is not necessarily a dense matrix. In fact, in the
case of ACOPF problems, $M_\text{cond}$ is still highly sparse, since
$W$ and $A$ are graph-induced banded systems.  Thus, exploiting
sparsity is still necessary to enable scalable computations.  In
general NLPs, however, the condensation strategy can arbitrarily
increase the density of the KKT system. Thus, the condensed-space IPM
strategy needs to be used with caution.

The reason that the condensation strategy is particularly relevant for
GPUs is that the matrix in \eqref{eqn:very-long-condensed} is positive
definite upon application of the standard inertia correction
method. Typically, to guarantee that the computed step direction is a
descent direction, we need a condition that
$\text{inertia}(M_\text{aug}) = (n+5m,0,m)$. Here, inertia refers to
the tuple of positive, zero, and negative eigenvalues. Accordingly, we
employ inertia correction methods to modify the augmented KKT system
so that the desired conditions on the inertia are satisfied.

By Sylvester's law, we have
\begin{align*}
  &\text{inertia}(M_\text{aug}) = (n+5m,0,m)\\
  &\iff \text{inertia}(M_\text{cond}) = (n,0,0).
\end{align*}
Thus, any choice of $\delta_w,\delta_c>0$ that makes the condensed KKT
system positive definite yields the desired inertia condition.
An important observation here is that the condensed KKT matrix
with desired inertia condition is always positive definite.
Thus, $M_{\text{cond}}$ can be factorized with fixed pivoting (e.g.,
Cholesky factorization or LU refactorization), which is significantly
more amenable to parallel implementation than is indefinite LBL$^\top$
factorization (standard in the state-of-the-art IPM
algorithms but requires the use of pivoting).

As we approach the solution, multiple constraints become active:
in the diagonal matrices $\Sigma_x$ and $\Sigma_s$, the terms associated
with the active (resp. inactive) variables $(x, s)$ go to infinity
(resp. $0$).
As such, the presence of active constraints can arbitrarily
increase the conditioning of the KKT system, leading to
an ill-conditioned KKT system in \eqref{eqn:very-long-condensed}.
Consequently, a single triangular solve may not provide a
sufficiently accurate step direction. Accordingly, iterative
refinement methods are employed to refine the solution by performing
multiple triangular solves. Notably, iterative refinement is applied
to the full KKT system, rather than the condensed system in
\eqref{eqn:very-long-condensed}.

\subsection{Line Search and IPM Iterations}

The step size can be determined by using the line search procedure.
Although numerous alternative approaches exist, we follow the
filter line search method implemented in the Ipopt solver. The line
search procedure employed here determines the step size by performing
a backtracking line search until a trial point satisfying sufficient
progress conditions is satisfied and acceptable by the
filter. Furthermore, in order to enhance the convergence behavior, various
additional strategies are implemented, such as the second-order
correction, restoration phase, and automatic scaling. For the details
of the implementation of the filter line search and various additional
strategies, the readers are referred to
\cite{wachter2006implementation}.

The step size and direction obtained above can be implemented as follows:
\begin{align}
  (x,s,y) &\leftarrow (x,s,y)+ \alpha (\Delta x, \Delta s, \Delta y),\label{eqn:iter}\\\nonumber
  (z_x^\flat, z_x^\sharp, z_s^\flat, z_s^\sharp) &\leftarrow (z_x^\flat, z_x^\sharp, z_s^\flat, z_s^\sharp) + \alpha_z (\Delta z_x^\flat, \Delta z_x^\sharp, \Delta z_s^\flat, \Delta z_s^\sharp).
\end{align}
The iteration in \eqref{eqn:iter} is repeated until the convergence
criterion is satisfied. The convergence criterion is defined as
\begin{align}\label{eqn:criteria}
\text{residual}(x^{(\ell)}, s^{(\ell)}, y^{(\ell)},
z^{(\ell)\flat}_x, z^{(\ell)\sharp}_x, z^{(\ell)\flat}_s,
z^{(\ell)\sharp}_s) <\epsilon_{\text{tol}},
\end{align}
where
$\text{residual}(\cdot)$ is a scaled version of the residual to the
first-order conditions in \eqref{eqn:first}. We summarize our
condensed-space IPM in Algorithm \ref{alg:con-ipm}.

\begin{algorithm}[t]
  \caption{Condensed-Space IPM}
  \label{alg:con-ipm}
  \begin{algorithmic}[1]
    \REQUIRE Primal-dual solution guesses $x,y, z^\flat, z^\sharp$, bounds $x^\flat,x^\sharp, s^\flat, s^\sharp$, callbacks $f(\cdot)$, $g(\cdot)$, $\nabla_x f(\cdot)$, $\nabla_x g(\cdot)$, $\nabla^2_{xx} \mathcal{L}(\cdot)$, and tolerance $\epsilon_{\text{tol}}$
    \STATE Relax the equality constraints by \eqref{eqn:relax} and initialize the slack $s$ and the associated dual variables $z^\flat_s, z^\sharp_s$.
    \WHILE{convergence criteria \eqref{eqn:criteria} not satisfied}
    \STATE Solve the condensed KKT system \eqref{eqn:very-long-condensed} with $\delta_w=\delta_c=0$ to compute the primal step $\Delta x$ and recover the dual steps $\Delta y, \Delta z_x^\flat, \Delta z_x^\sharp, \Delta z_s^\flat, \Delta z_s^\sharp$ by \eqref{eqn:recover-1} and \eqref{eqn:recover-2}.
    \STATE Determine the need for regularization and, if necessary, recompute the step directions with proper choices of $\delta_w,\delta_c>0$.
    \STATE Choose step sizes $\alpha,\alpha_z>0$ via line search.
    \STATE Update the solution by \eqref{eqn:iter}.
    \STATE Update filter and barrier parameter $\mu$.
    \ENDWHILE
    \RETURN The first-order stationary points $x^\star,y^\star, z^{\flat\star}, z^{\sharp\star}$
  \end{algorithmic}
\end{algorithm}

\subsection{Notes on the Implementation}
We have implemented the condensed-space IPM by adapting our code base
in MadNLP.jl, a port of Ipopt in Julia.
A key feature of MadNLP is that the IPM is implemented with a
high level of abstraction, while the specific handling of the data
structures within the KKT systems is carried out by data-type specific
kernel functions. This design allows us to apply the mathematical
operations equivalent to Ipopt to different KKT data structures, such as
{\tt SparseKKTSystem}, {\tt DenseKKTSystem}, and {\tt
DenseCondensedKKTSystem},  whose data are stored either on the host or in device memory. For the implementation of the condensed-space IPM
presented in this paper, we have added a new type of KKT system called
{\tt SparseCondensedKKTSystem} and implemented additional kernels
needed for handling the data structures specific to this KKT system
type. This approach ensures that we are performing mathematically
equivalent operations as in the mature, extensively tested existing
code base. This also allows us to easily switch between different KKT
system types, which is crucial for experimenting with various solvers
and data structures, as well as for efficiently leveraging GPU
acceleration when available. Furthermore, by maintaining this level of
abstraction, the condensed-space IPM can be seamlessly integrated into
the existing framework, making it easier to maintain and extend in the
future.

\section{Numerical Results}\label{sec:num}

This section presents the numerical benchmark results, comparing our
method against state-of-the-art methods on CPUs for solving standard
ACOPF problems.

\subsection{Methods}

We compared four different configurations of NLP solution frameworks:
\begin{align}
  \label{config-1}\tag{Config 1} \bullet\;&\text{MadNLP.jl+ExaModels.jl+cuSOLVER (GPU)}\\
  \label{config-2}\tag{Config 2} \bullet\;&\text{MadNLP.jl+ExaModels.jl+Ma27 (CPU)}\\
  \label{config-3}\tag{Config 3} \bullet\;&\text{Ipopt+AMPL+Ma27 (CPU)}\\
  \label{config-4}\tag{Config 4} \bullet\;&\text{Ipopt+JuMP.jl+Ma27 (CPU)}.
\end{align}

\ref{config-1} is our main GPU configuration; \ref{config-2}
represents our implementation running on CPU; and \ref{config-3} and
\ref{config-4} are used as benchmarks. \ref{config-1} and
\ref{config-2} share a significant amount of code, especially the
high-level abstractions, but they differ in how they handle the KKT
systems. In \ref{config-1}, MadNLP.jl applies the condensed-space IPM
along with the inequality relaxation strategy, while in
\ref{config-2}, MadNLP.jl applies IPM based on the indefinite,
noncondensed KKT system, as in \eqref{eqn:kkt-indefinite}. In
\ref{config-1}, we use the cuSOLVER library to solve the condensed KKT
system via Cholesky factorization. The numerical
factorization and triangular solves are performed by cuSOLVER with the
fixed pivot sequence obtained with an approximate minimum degree ordering algorithm \cite{amestoyApproximateMinimumDegree1996}, implemented in AMD.jl \cite{montoisonAMDJlJulia2020}.  Software and hardware details of each
configuration are illustrated in Table \ref{tbl:settings}. The ACOPF
problem is formulated by using the model from the rosetta-opf project
\cite{rosetta-opf}, which is based on the models in PowerModels.jl \cite{8442948}, and the test cases are obtained from the pglib-opf
repository \cite{babaeinejadsarookolaee2019power}. We have selected
the goc and pegase cases because they contain large-scale instances.  The
external packages are called from Julia, through thin wrapper
packages, such as Ipopt.jl and AmplNLWriters.jl. A tolerance of
$10^{-4}$ is set for MadNLP.jl and Ipopt solvers, with other solver
options adjusted to ensure a fair comparison across different
solvers. The results can be reproduced with the script available at
\url{https://github.com/sshin23/opf-on-gpu}.

\subsection{Results}

\begin{figure}[t]
  \centering
  \includegraphics[width=.4\textwidth]{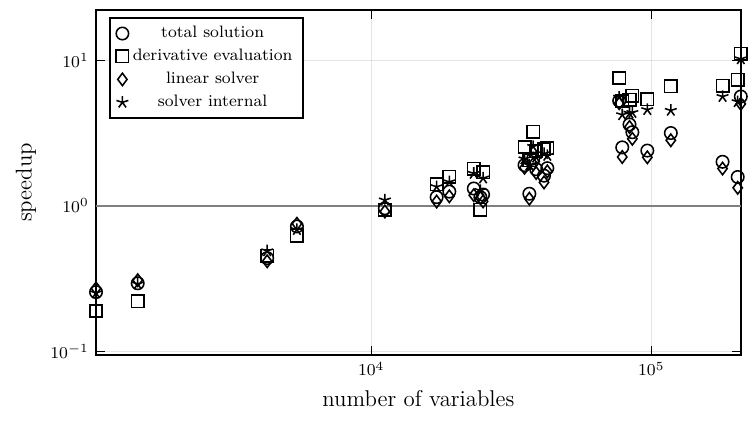}
  \vspace{-.1in}
  \caption{Speedup achieved by using GPUs.}
  \label{fig:speedup}
  \vspace{-.2in}
\end{figure}

The numerical benchmark results, including total solution time and its
breakdown into linear algebra and derivative evaluation time (with the
remainder considered as solver internal time), are shown in Table
\ref{tbl:results}. The quality of the solution (objective value and
constraint violation measured by $\|\cdot\|_\infty$) is shown in Table
\ref{tbl:quality}. Figure \ref{fig:speedup} visually represents the
speedup brought by GPUs, by comparing the timing results of
\ref{config-1} and \ref{config-2}.

\begin{table*}[t]
  \scriptsize
  \centering
  \caption{Details of Numerical Experiment Settings}
  \begin{tabular}{|l|c|c|c|c|}
    \hline
    & {\textbf{MadNLP+ExaModels+cuSOLVER}}
    & {\textbf{MadNLP+ExaModels+Ma27}}
    & {\textbf{Ipopt+AMPL+Ma27}}
    & {\textbf{Ipopt+JuMP+Ma27}}\\
    &\textbf{(GPU)} &\textbf{(CPU)} &\textbf{(CPU)}& \textbf{(CPU)}\\
    \hline
    \textbf{Optimization Solver} & \multicolumn{2}{c|}{MadNLP.jl (dev)$^*$} & \multicolumn{2}{c|}{Ipopt (v3.13.3)} \\
    \hline
    \textbf{Derivative Evaluations} & \multicolumn{2}{c|}{ExaModels.jl (dev)$^*$} &  AMPL Solver Library & JuMP.jl (v1.12.0)\\
    \hline
    \textbf{Linear Solver} &  cuSOLVER (v11.4.5) &\multicolumn{3}{c|}{Ma27 (v2015.06.23)}\\
    \hline
    \textbf{Hardware} & NVIDIA Quadro GV100 & \multicolumn{3}{c|}{Intel Xeon Gold 6140}\\
    \hline
  \end{tabular}\\
  $^*$Specific commit hashes are available at \url{https://github.com/sshin23/opf-on-gpu}
  \label{tbl:settings}
  \vspace{-.1in}
\end{table*}
\begin{table*}[t]
  \scriptsize
  \centering
  \caption{Numerical Results}
  \begin{tabular}{|l|c|c|cccc|cccc|ccc|ccc|}
  \hline
  \multirow{3}{*}{\textbf{Case}}
  & \multirow{3}{*}{nvars}
  & \multirow{3}{*}{ncons}
  & \multicolumn{4}{c|}{\textbf{MadNLP+ExaModels+cuSOLVER}}
  & \multicolumn{4}{c|}{\textbf{MadNLP+ExaModels+Ma27}}
  & \multicolumn{3}{c|}{\textbf{Ipopt+AMPL+Ma27}}
  & \multicolumn{3}{c|}{\textbf{Ipopt+JuMP+Ma27}}\\
  & & &\multicolumn{4}{c|}{\textbf{(GPU)}} &\multicolumn{4}{c|}{\textbf{(CPU)}} &\multicolumn{3}{c|}{\textbf{(CPU)}}&\multicolumn{3}{c|}{\textbf{(CPU)}}
  \\
  \cline{4-17}
  & & 
  & iter & deriv.$^\dag$ & lin.$^\dag$ & total$^\dag$
  & iter & deriv.$^\dag$ & lin.$^\dag$ & total$^\dag$
  & iter & deriv.$^\ddag$ & total$^\ddag$
  & iter & deriv.$^\ddag$ & total$^\ddag$
  \\
  \hline
89\_pegase 
&   1.0k
&   1.6k
& 28 
&  0.02
&  0.12
&  0.22
& 30 
&  0.00
&  0.03
& \bf{ 0.06}
& 29 
&  0.04
&  0.09
& 29 
&  0.12
&  0.18
\\

179\_goc 
&   1.5k
&   2.2k
& 30 
&  0.03
&  0.17
&  0.30
& 43 
&  0.01
&  0.05
& \bf{ 0.09}
& 42 
&  0.05
&  0.13
& 42 
&  0.17
&  0.26
\\

500\_goc 
&   4.3k
&   6.1k
& 36 
&  0.04
&  0.31
&  0.47
& 35 
&  0.02
&  0.13
& \bf{ 0.20}
& 36 
&  0.14
&  0.31
& 34 
&  0.43
&  0.64
\\

793\_goc 
&   5.4k
&   8.0k
& 33 
&  0.03
&  0.21
&  0.33
& 31 
&  0.02
&  0.16
& \bf{ 0.24}
& 31 
&  0.20
&  0.39
& 30 
&  0.58
&  0.82
\\

1354\_pegase 
&  11.2k
&  16.6k
& 44 
&  0.06
&  0.48
&  0.73
& 45 
&  0.06
&  0.44
&  {\bf 0.70}
& 41 
&  0.94
&  1.48
& 41 
&  2.40
&  3.04
\\
\hline
2312\_goc 
&  17.1k
&  25.7k
& 38 
&  0.06
&  0.75
& \bf{ 1.02}
& 40 
&  0.08
&  0.80
&  1.16
& 38 
&  1.45
&  2.33
& 38 
&  3.04
&  4.05
\\

2000\_goc 
&  19.0k
&  29.4k
& 36 
&  0.06
&  0.76
& \bf{ 1.05}
& 38 
&  0.09
&  0.88
&  1.32
& 39 
&  1.72
&  2.79
& 38 
&  5.20
&  6.41
\\

3022\_goc 
&  23.2k
&  35.0k
& 43 
&  0.08
&  1.09
& \bf{ 1.47}
& 49 
&  0.14
&  1.29
&  1.93
& 47 
&  2.57
&  4.02
& 47 
&  7.49
&  9.16
\\

2742\_goc 
&  24.5k
&  38.2k
& 151 
&  0.50
&  4.97
& \bf{ 6.67}
& 122 
&  0.46
&  5.63
&  7.63
& 98 
&  8.50
& 13.91
& 99 
& 21.23
& 27.04
\\

2869\_pegase 
&  25.1k
&  37.8k
& 52 
&  0.10
&  1.44
& \bf{ 1.90}
& 52 
&  0.16
&  1.54
&  2.27
& 50 
&  3.27
&  4.99
& 50 
&  6.24
&  8.10
\\
\hline
3970\_goc 
&  35.3k
&  54.4k
& 44 
&  0.09
&  1.62
& \bf{ 2.05}
& 45 
&  0.22
&  2.95
&  3.90
& 60 
&  5.42
&  9.94
& 43 
&  7.36
& 10.95
\\

4020\_goc 
&  36.7k
&  57.0k
& 70 
&  0.14
&  5.50
& \bf{ 6.19}
& 59 
&  0.30
&  6.14
&  7.48
& 55 
&  5.28
& 11.66
& 55 
& 10.75
& 17.33
\\

4917\_goc 
&  37.9k
&  56.9k
& 48 
&  0.09
&  1.47
& \bf{ 1.93}
& 57 
&  0.28
&  2.83
&  4.07
& 53 
&  5.12
&  7.98
& 53 
&  9.80
& 13.04
\\

4601\_goc 
&  38.8k
&  59.6k
& 71 
&  0.15
&  2.77
& \bf{ 3.46}
& 66 
&  0.35
&  4.66
&  6.17
& 69 
&  7.02
& 12.72
& 68 
& 13.37
& 19.12
\\

4837\_goc 
&  41.4k
&  64.0k
& 57 
&  0.13
&  2.69
& \bf{ 3.31}
& 56 
&  0.32
&  3.89
&  5.32
& 56 
&  8.22
& 13.09
& 56 
& 12.36
& 17.13
\\
\hline
4619\_goc 
&  42.5k
&  66.3k
& 54 
&  0.11
&  2.84
& \bf{ 3.40}
& 46 
&  0.27
&  4.89
&  6.15
& 48 
&  8.30
& 14.14
& 46 
& 10.37
& 15.57
\\

10000\_goc 
&  76.8k
& 112.4k
& 56 
&  0.10
&  1.93
& \bf{ 2.53}
& 77 
&  0.76
&  9.81
& 13.30
& 74 
& 14.56
& 24.85
& 74 
& 24.71
& 36.00
\\

8387\_pegase 
&  78.7k
& 118.7k
& 67 
&  0.14
&  4.25
& \bf{ 5.05}
& 70 
&  0.75
&  9.19
& 12.72
& 69 
& 16.70
& 26.54
& 69 
& 25.97
& 36.19
\\

9591\_goc 
&  83.6k
& 130.6k
& 69 
&  0.15
&  5.21
& \bf{ 6.03}
& 65 
&  0.78
& 18.12
& 21.81
& 64 
& 16.92
& 38.50
& 62 
& 34.96
& 54.47
\\

9241\_pegase 
&  85.6k
& 130.8k
& 63 
&  0.13
&  3.38
&  {\bf 4.16}
& 63 
&  0.74
&  9.76
& 13.31
& 61 
& 15.87
& 26.66
& 61 
& 25.41
& 36.69
\\
\hline
10480\_goc 
&  96.8k
& 150.9k
& 70 
&  0.17
&  8.88
&  {\bf 9.80}
& 66 
&  0.90
& 19.11
& 23.46
& 64 
& 17.58
& 38.82
& 63 
& 31.76
& 52.65
\\

13659\_pegase 
& 117.4k
& 170.6k
& 66 
&  0.14
&  4.46
&  {\bf 5.36}
& 58 
&  0.92
& 12.56
& 16.94
& 64 
& 19.90
& 35.79
& 64 
& 35.97
& 53.02
\\

19402\_goc 
& 179.6k
& 281.7k
& 102 
&  0.29
& 30.46
& {\bf 32.08}
& 70 
&  1.93
& 54.88
& 64.29
& 70 
& 36.34
& 94.72
& 70 
& 65.25
& 121.72
\\

24464\_goc 
& 203.4k
& 313.6k
& 80 
&  0.25
& 25.11
& {\bf 26.69}
& 58 
&  1.81
& 33.33
& 42.03
& 58 
& 34.33
& 71.25
& 58 
& 61.04
& 99.47
\\

30000\_goc 
& 208.6k
& 307.8k
& 153 
&  0.43
& 14.86
& {\bf 16.79}
& 136 
&  4.74
& 74.54
& 94.62
& 180 
& 105.03
& 248.64
& 126 
& 133.13
& 206.70

  \\
  \hline
\end{tabular}\\
  $^\dag$Wall time (sec) measured by Julia. $^\ddag$CPU time (sec) reported by Ipopt.
  \label{tbl:results}
\end{table*}
\begin{table*}[t]
  \scriptsize
  \centering
  \caption{Solution Quality}
  \label{tbl:quality}
  \begin{tabular}{|l|cc|cc|cc|cc|}
  \hline
  \multirow{3}{*}{\textbf{Case}}
  & \multicolumn{2}{c|}{\textbf{MadNLP+ExaModels+cuSOLVER}}
  & \multicolumn{2}{c|}{\textbf{MadNLP+ExaModels+Ma27}}
  & \multicolumn{2}{c|}{\textbf{Ipopt+AMPL+Ma27}}
  & \multicolumn{2}{c|}{\textbf{Ipopt+JuMP+Ma27}}\\
  &\multicolumn{2}{c|}{\textbf{(GPU)}} &\multicolumn{2}{c|}{\textbf{(CPU)}} &\multicolumn{2}{c|}{\textbf{(CPU)}}&\multicolumn{2}{c|}{\textbf{(CPU)}}
  \\
  \cline{2-9}
  & objective & constr. viol.
  & objective & constr. viol.
  & objective & constr. viol.
  & objective & constr. viol.
  \\
  \hline
89\_pegase 
& 1.07023029e+05
& 1.69977362e-03
& 1.07277300e+05
& 1.69995406e-03
& 1.07273132e+05
& 1.69762454e-02
& 1.07273132e+05
& 1.69762454e-02
\\

179\_goc 
& 7.54098231e+05
& 3.64045772e-03
& 7.54215279e+05
& 3.64095371e-03
& 7.54214091e+05
& 1.05727439e-02
& 7.54214091e+05
& 1.05727439e-02
\\

500\_goc 
& 4.53056588e+05
& 1.16442922e-03
& 4.54894607e+05
& 1.16461929e-03
& 4.54894301e+05
& 1.16449188e-03
& 4.54894349e+05
& 1.16443248e-03
\\

793\_goc 
& 2.59660004e+05
& 1.12495280e-03
& 2.60179408e+05
& 1.14373500e-03
& 2.60177953e+05
& 2.52890328e-02
& 2.60177960e+05
& 2.52825510e-02
\\

1354\_pegase 
& 1.25574315e+06
& 4.18838427e-03
& 1.25874608e+06
& 4.18894441e-03
& 1.25873160e+06
& 2.91106529e-02
& 1.25873160e+06
& 2.91106529e-02
\\
\hline
2312\_goc 
& 4.40492687e+05
& 1.95782217e-03
& 4.41301927e+05
& 1.98487972e-03
& 4.41301012e+05
& 2.86441953e-03
& 4.41301012e+05
& 2.86441953e-03
\\

2000\_goc 
& 9.66186544e+05
& 1.07957382e-03
& 9.73392385e+05
& 1.07991565e-03
& 9.73392524e+05
& 1.07970410e-03
& 9.73392602e+05
& 1.07958552e-03
\\

3022\_goc 
& 6.00461469e+05
& 1.60590210e-03
& 6.01341340e+05
& 1.92264271e-03
& 6.01340934e+05
& 7.06720510e-03
& 6.01340934e+05
& 7.06720510e-03
\\

2742\_goc 
& 2.70328757e+05
& 9.99725733e-04
& 2.75672815e+05
& 9.99997332e-04
& 2.75672759e+05
& 1.13868333e-03
& 2.75672759e+05
& 1.13868340e-03
\\

2869\_pegase 
& 2.45584120e+06
& 4.18833905e-03
& 2.46259584e+06
& 4.18882610e-03
& 2.46258759e+06
& 3.15283321e-02
& 2.46258759e+06
& 3.15283321e-02
\\
\hline
3970\_goc 
& 9.27998953e+05
& 6.41922608e-04
& 9.60666837e+05
& 6.42469892e-04
& 9.60667021e+05
& 6.42371530e-04
& 9.60667776e+05
& 6.41960999e-04
\\

4020\_goc 
& 8.02565861e+05
& 1.29969745e-03
& 8.21952202e+05
& 1.29999868e-03
& 8.21952543e+05
& 1.29986624e-03
& 8.21952543e+05
& 1.29986624e-03
\\

4917\_goc 
& 1.38537252e+06
& 1.54172485e-03
& 1.38769645e+06
& 1.70860688e-03
& 1.38769342e+06
& 1.62739725e-02
& 1.38769342e+06
& 1.62739725e-02
\\

4601\_goc 
& 7.92510931e+05
& 9.99886244e-04
& 8.25898288e+05
& 9.99978318e-04
& 8.25898470e+05
& 9.99896654e-04
& 8.25898481e+05
& 9.99894295e-04
\\

4837\_goc 
& 8.60071647e+05
& 9.92673673e-04
& 8.72192598e+05
& 9.92934504e-04
& 8.72192733e+05
& 9.92677263e-04
& 8.72192733e+05
& 9.92677263e-04
\\
\hline
4619\_goc 
& 4.66738422e+05
& 8.80364611e-04
& 4.76659294e+05
& 8.80485073e-04
& 4.76659432e+05
& 8.80367536e-04
& 4.76659432e+05
& 8.80367536e-04
\\

10000\_goc 
& 1.34739992e+06
& 5.36209615e-04
& 1.35370965e+06
& 5.40993748e-04
& 1.35371078e+06
& 6.56672045e-04
& 1.35371173e+06
& 6.56367359e-04
\\

8387\_pegase 
& 2.74980910e+06
& 9.99884691e-03
& 2.77083829e+06
& 9.99896893e-03
& 2.77062704e+06
& 5.30460965e-02
& 2.77062704e+06
& 5.30460965e-02
\\

9591\_goc 
& 1.02516095e+06
& 9.91659468e-04
& 1.06148769e+06
& 9.91997903e-04
& 1.06148806e+06
& 9.91795084e-04
& 1.06148807e+06
& 9.91788322e-04
\\

9241\_pegase 
& 6.21773526e+06
& 4.18380647e-03
& 6.24208171e+06
& 4.18787958e-03
& 6.24207325e+06
& 3.76440386e-02
& 6.24207325e+06
& 3.76440386e-02
\\
\hline
10480\_goc 
& 2.27696973e+06
& 1.09983709e-03
& 2.31442783e+06
& 1.09996886e-03
& 2.31442450e+06
& 1.67932256e-02
& 2.31442450e+06
& 1.67932256e-02
\\

13659\_pegase 
& 8.92385389e+06
& 1.99904428e-03
& 8.94679835e+06
& 1.99980680e-03
& 8.94680070e+06
& 1.54477837e-02
& 8.94680070e+06
& 1.54477837e-02
\\

19402\_goc 
& 1.93394723e+06
& 1.19983797e-03
& 1.97755237e+06
& 1.19999867e-03
& 1.97755235e+06
& 1.19986568e-03
& 1.97755235e+06
& 1.19986568e-03
\\

24464\_goc 
& 2.58935629e+06
& 7.24722104e-04
& 2.62932336e+06
& 7.24944021e-04
& 2.62932439e+06
& 7.24724162e-04
& 2.62932439e+06
& 7.24724162e-04
\\

30000\_goc 
& 1.11353160e+06
& 1.40161701e-03
& 1.14190983e+06
& 1.40292333e-03
& 1.14191122e+06
& 1.40225897e-03
& 1.14190714e+06
& 1.40184075e-03

  \\
  \hline
\end{tabular}\\
\end{table*}

\paragraph*{Convergence Pattern}
When comparing the solvers' performance in terms of the
IPM iteration counts, MadNLP.jl is as
efficient as the state-of-the-art solver Ipopt. The IPM iteration
count is nearly the same as that of Ipopt for achieving the same level
of accuracy in the final solution (see Table \ref{tbl:quality}).
{\color{black} Note that
the constraint violation is not strictly less than the tolerance in all configurations
as relative (scaled) constraint violation is used for the convergence criteria.}
This
suggests that running mathematically equivalent operations on GPUs (by
using a shared code base in high-level abstractions) can yield a similar
degree of effectiveness in terms of IPM convergence.

\paragraph*{Performance of AD}
Next, we discuss the effectiveness of parallel AD on GPUs. We observe
that even on CPUs, ExaModels.jl is substantially faster than the AD
routines implemented in AMPL or JuMP.jl. Indeed, ExaModels.jl generates
derivative functions specifically compiled for the type of model,
including optimization for the distinctive computational pattern found
in the model. When comparing ExaModels.jl running on CPUs and GPUs, we
observe a further speedup on the GPU of up to 10x for large
instances (e.g., case30000\_goc); remarkably, this is 300 times faster
than JuMP.jl. This demonstrates that adopting SIMD abstraction and
parallelizing AD brings significant computational gain. 

\paragraph*{Performance of the Condensed-Space IPM}
We next discuss the linear solver time. While the speedup achieved by
linear solvers is only moderate, this has a high impact on the overall
speedup because the linear solver time constitutes a significant
portion of the total solution time. {\color{black} Figure
\ref{fig:speedup} reveals that linear algebra is the only
computational bottleneck for large-scale instances. While derivative
evaluation and solver internal computation could achieve a 10x speedup,
the linear algebra part can only achieve a 5x speedup even for the
largest case. For case24464\_goc, the GPU linear solver performance
was close to that of the CPU solver.  The investigation of under which
circumstances cuSOLVER is more effective warrants further research.}

Solver internal time could also be significantly accelerated through
GPU utilization. We can observe that the speedup in solver internal
operations is consistently greater than the speedup in linear
solvers. Because of the frequent use of \ref{eqn:pattern-2} and
\ref{eqn:pattern-3} operations,  however, the speedup in solver internal
operations is less than that of derivative evaluations.

Overall, our GPU implementation exhibits significant speedup across
all components (derivative evaluation, linear algebra, and solver
internal computation) resulting in substantial gains in total solution
time. The results indicate that GPUs become more effective for
large-scale instances, particularly when the number of variables is
greater than 20,000. {\color{black} This is because the benefit of parallelization
is greater for operations over large data.}
Notably, for the largest instance,
case30000\_goc, our GPU implementation is 4 times faster than our CPU
implementation and approximately 10 times faster than
state-of-the-art tools (Ipopt, JuMP.jl, and Ma27).
This demonstrates that GPUs can bring significant computational
gains for large-scale AC OPF problems, enabling the solution of
previously inconceivable problems due to the limitations of CPU-based
solution tools.

\paragraph*{Known Limitations}
The key limitation of our method is the decreased precision of the
final solution. Reliable convergence is achieved only up to a
tolerance of $10^{-4}$. We have observed that the condensation of the
KKT system worsens the conditioning of the already ill-conditioned
augmented KKT system, resulting in even higher condition numbers,
particularly when the solution is almost converged. For instance, in
the case of case118\_ieee, the condition number of the KKT system at
the solution (with tolerance set to $10^{-4}$) is $9.43\times 10^{11}$
for the augmented KKT system and $3.15\times 10^{14}$ for the
condensed KKT system. Consequently, the achievable precision of
the solution is reduced. Further investigation into the possibility of
obtaining higher precision will be necessary.

\section{Conclusions and Future Outlook}\label{sec:conc}
We have presented an NLP solution framework for solving large-scale AC
OPF problems. By leveraging the SIMD abstraction of NLPs and a
condensed-space IPM, we have effectively eliminated the need for
serial computations, enabling the implementation of a solution
framework that can run entirely on GPUs. Our method has demonstrated
promising results, achieving a 5x speedup when compared with CPU
implementations for large-scale ACOPF problems. Notably, our approach
outperforms one of the state-of-the-art CPU-based implementations by
a factor of 12. These results, along with our packages
MadNLP.jl and ExaModels.jl, showcase a significant advancement in our
capabilities in dealing with large-scale optimization problems in
power systems and underscore the potential of accelerated computing in
large-scale optimization area. However, the condensation procedure
leads to an increase in the condition number of the KKT system,
resulting in decreased final solution accuracy. Addressing the
challenges posed by ill-conditioning remains an important aspect
of future work. In the following paragraphs, we discuss some
remaining open questions and future outlooks.

\paragraph*{Obtaining Higher Numerical Precision}
While we have focused on the IPM, other constrained optimization
paradigms, such as penalty methods and augmented Lagrangian methods,
exist, and similar strategies based on condensed linear systems can be
developed. It would be valuable to investigate which algorithm would
be the right paradigm for constrained large-scale optimization on GPUs
that can best handle the ill-conditioning issue of the condensed KKT
system and, in turn, achieve the highest degree of accuracy.

\paragraph*{Security-Constrained, Multiperiod, Distribution OPFs}
Although the proposed method has demonstrated significant computational
advantages for transmission ACOPF problems, our results can also
be interpreted that efficient CPUs can still handle these problems
reasonably well. We anticipate that there will be more substantial
performance gains for larger-scale optimization problems, such as
security-constrained and multiperiod OPFs or joint optimization
problems involving transmission, distribution, and gas network
systems. We are interested in exploring Schur complement-based
decomposition approaches, combined with the condensation-based
strategy, similarly to \cite{pacaud2023parallel}, to demonstrate even
greater scalability.

\paragraph*{Alternative Linear Solvers}
While cuSOLVER has been effective for solving the condensed KKT
systems using LU factorization, Cholesky factorization holds promise
for better performance due to lower computational complexity and the
ability to reveal the inertia of the KKT system. We are
interested in exploring other linear solver options, such as CHOLMOD
\cite{chen2008algorithm}, Baspacho \cite{pineda2022theseus}, and HyKKT
\cite{regev2023hykkt}.

\paragraph*{Portability}
Our implementation is currently  tested only on NVIDIA GPUs, but our
GPU implementation is largely based on array programming and
KernelAbstractions.jl in Julia, which are in principle  compatible with
various GPU architectures, including AMD, Intel, and Apple GPUs. By
incorporating cross-architecture linear solvers, we envision
supporting a broader class of GPU accelerators in the future.

\section*{Acknowledgment}
This material is based upon work supported by the U.S. Department of Energy, Office of Science, Office of Advanced Scientific Computing Research (ASCR) under Contract DE-AC02-06CH11347.

\bibliographystyle{ieeetr}
\bibliography{main}

\begin{thebibliography}{10}

\bibitem{anitescu2021targeting}
M.~Anitescu, K.~Kim, Y.~Kim, A.~Maldonado, F.~Pacaud, V.~Rao, M.~Schanen,
  S.~Shin, and A.~Subramanian, ``Targeting {E}xascale with {J}ulia on
  {G}{P}{U}s for multiperiod optimization with scenario constraints,'' {\em
  SIAG/OPT Views and News}, 2021.

\bibitem{simdiff}
S.~Shin, ``{E}xa{M}odels.jl.'' \url{https://github.com/sshin23/ExaModels.jl}.

\bibitem{madnlp}
S.~Shin and F.~Pacaud, ``{M}ad{NLP}.jl.''
  \url{https://github.com/MadNLP/MadNLP.jl}.

\bibitem{dunning2017jump}
I.~Dunning, J.~Huchette, and M.~Lubin, ``{JuMP}: A modeling language for
  mathematical optimization,'' {\em SIAM review}, vol.~59, no.~2, pp.~295--320,
  2017.

\bibitem{fourer1990modeling}
R.~Fourer, D.~M. Gay, and B.~W. Kernighan, ``A modeling language for
  mathematical programming,'' {\em Management Science}, vol.~36, no.~5,
  pp.~519--554, 1990.

\bibitem{swirydowicz2022linear}
K.~{\'S}wirydowicz, E.~Darve, W.~Jones, J.~Maack, S.~Regev, M.~A. Saunders,
  S.~J. Thomas, and S.~Pele{\v{s}}, ``Linear solvers for power grid
  optimization problems: a review of {GPU}-accelerated linear solvers,'' {\em
  Parallel Computing}, vol.~111, p.~102870, 2022.

\bibitem{nocedal2006numerical}
J.~Nocedal and S.~J. Wright, {\em Numerical optimization}.
\newblock Springer, 2006.

\bibitem{babaeinejadsarookolaee2019power}
S.~Babaeinejadsarookolaee, A.~Birchfield, R.~D. Christie, C.~Coffrin,
  C.~DeMarco, R.~Diao, M.~Ferris, S.~Fliscounakis, S.~Greene, R.~Huang, {\em
  et~al.}, ``The power grid library for benchmarking {AC} optimal power flow
  algorithms,'' {\em arXiv preprint arXiv:1908.02788}, 2019.

\bibitem{cao2016augmented}
Y.~Cao, A.~Seth, and C.~D. Laird, ``An augmented {Lagrangian} interior-point
  approach for large-scale {NLP} problems on graphics processing units,'' {\em
  Computers \& Chemical Engineering}, vol.~85, pp.~76--83, 2016.

\bibitem{pacaud2023parallel}
F.~Pacaud, M.~Schanen, S.~Shin, D.~A. Maldonado, and M.~Anitescu, ``Parallel
  interior-point solver for block-structured nonlinear programs on {SIMD/GPU}
  architectures,'' {\em arXiv preprint arXiv:2301.04869}, 2023.

\bibitem{pacaud2022feasible}
F.~Pacaud, D.~A. Maldonado, S.~Shin, M.~Schanen, and M.~Anitescu, ``A feasible
  reduced space method for real-time optimal power flow,'' {\em Electric Power
  Systems Research}, vol.~212, p.~108268, 2022.

\bibitem{pacaud2023accelerating}
F.~Pacaud, S.~Shin, M.~Schanen, D.~A. Maldonado, and M.~Anitescu,
  ``Accelerating condensed interior-point methods on {SIMD/GPU}
  architectures,'' {\em Journal of Optimization Theory and Applications},
  pp.~1--20, 2023.

\bibitem{lin1999newton}
C.-J. Lin and J.~J. Mor{\'e}, ``Newton's method for large bound-constrained
  optimization problems,'' {\em SIAM Journal on Optimization}, vol.~9, no.~4,
  pp.~1100--1127, 1999.

\bibitem{kim2022accelerated}
Y.~Kim and K.~Kim, ``Accelerated computation and tracking of {AC} optimal power
  flow solutions using {GPU}s,'' in {\em Workshop Proceedings of the 51st
  International Conference on Parallel Processing}, pp.~1--8, 2022.

\bibitem{kim2021leveraging}
Y.~Kim, F.~Pacaud, K.~Kim, and M.~Anitescu, ``Leveraging {GPU} batching for
  scalable nonlinear programming through massive lagrangian decomposition,''
  {\em arXiv preprint arXiv:2106.14995}, 2021.

\bibitem{hiop_techrep}
C.~G. Petra, N.~Chiang, and J.~Wang, ``{HiOp} -- {U}ser {G}uide,'' Tech. Rep.
  LLNL-SM-743591, Center for Applied Scientific Computing, Lawrence Livermore
  National Laboratory, 2018.

\bibitem{regev2023hykkt}
S.~Regev, N.-Y. Chiang, E.~Darve, C.~G. Petra, M.~A. Saunders,
  K.~{\'S}wirydowicz, and S.~Pele{\v{s}}, ``{HyKKT}: a hybrid direct-iterative
  method for solving {KKT} linear systems,'' {\em Optimization Methods and
  Software}, vol.~38, no.~2, pp.~332--355, 2023.

\bibitem{Gravity}
H.~Hijazi, G.~Wang, and C.~Coffrin, ``Gravity: A mathematical modeling language
  for optimization and machine learning,'' {\em Machine Learning Open Source
  Software Workshop at NeurIPS 2018}, 2018.
\newblock Available at \url{www.gravityopt.com}.

\bibitem{jerez2012sparse}
J.~L. Jerez, E.~C. Kerrigan, and G.~A. Constantinides, ``A sparse and condensed
  {QP} formulation for predictive control of {LTI} systems,'' {\em Automatica},
  vol.~48, no.~5, pp.~999--1002, 2012.

\bibitem{cole2023exploiting}
D.~Cole, S.~Shin, F.~Pacaud, V.~M. Zavala, and M.~Anitescu, ``Exploiting
  {GPU/SIMD} architectures for solving linear-quadratic {MPC} problems,'' in
  {\em 2023 American Control Conference (ACC)}, pp.~3995--4000, IEEE, 2023.

\bibitem{chiang2014structured}
N.~Chiang, C.~G. Petra, and V.~M. Zavala, ``Structured nonconvex optimization
  of large-scale energy systems using {PIPS-NLP},'' in {\em 2014 Power Systems
  Computation Conference}, pp.~1--7, IEEE, 2014.

\bibitem{rodriguez2023scalable}
J.~S. Rodriguez, R.~B. Parker, C.~D. Laird, B.~L. Nicholson, J.~D. Siirola, and
  M.~L. Bynum, ``Scalable parallel nonlinear optimization with {PyNumero} and
  {P}arapint,'' {\em INFORMS Journal on Computing}, vol.~35, no.~2,
  pp.~509--517, 2023.

\bibitem{shin2021graph}
S.~Shin, C.~Coffrin, K.~Sundar, and V.~M. Zavala, ``Graph-based modeling and
  decomposition of energy infrastructures,'' {\em IFAC-PapersOnLine}, vol.~54,
  no.~3, pp.~693--698, 2021.

\bibitem{wachter2006implementation}
A.~W{\"a}chter and L.~T. Biegler, ``On the implementation of an interior-point
  filter line-search algorithm for large-scale nonlinear programming,'' {\em
  Mathematical Programming}, vol.~106, pp.~25--57, 2006.

\bibitem{bezanson2017julia}
J.~Bezanson, A.~Edelman, S.~Karpinski, and V.~B. Shah, ``Julia: A fresh
  approach to numerical computing,'' {\em SIAM Review}, vol.~59, no.~1,
  pp.~65--98, 2017.

\bibitem{amestoyApproximateMinimumDegree1996}
P.~R. Amestoy, T.~A. Davis, and I.~S. Duff, ``An {{Approximate Minimum Degree
  Ordering Algorithm}},'' {\em SIAM Journal on Matrix Analysis and
  Applications}, vol.~17, pp.~886--905, Oct. 1996.

\bibitem{montoisonAMDJlJulia2020}
A.~Montoison, D.~Orban, A.~S. Siqueira, and {contributors}, ``{{AMD}}.jl: {{A
  Julia}} interface to the {{AMD}} library of {{Amestoy}}, {{Davis}} and
  {{Duff}},'' May 2020.

\bibitem{rosetta-opf}
``rosetta-opf.'' \url{https://github.com/lanl-ansi/rosetta-opf}.

\bibitem{8442948}
C.~Coffrin, R.~Bent, K.~Sundar, Y.~Ng, and M.~Lubin, ``Power{M}odels.jl: An
  open-source framework for exploring power flow formulations,'' in {\em 2018
  Power Systems Computation Conference (PSCC)}, pp.~1--8, June 2018.

\bibitem{chen2008algorithm}
Y.~Chen, T.~A. Davis, W.~W. Hager, and S.~Rajamanickam, ``Algorithm 887:
  {CHOLMOD}, supernodal sparse {Cholesky} factorization and update/downdate,''
  {\em ACM Transactions on Mathematical Software (TOMS)}, vol.~35, no.~3,
  pp.~1--14, 2008.

\bibitem{pineda2022theseus}
L.~Pineda, T.~Fan, M.~Monge, S.~Venkataraman, P.~Sodhi, R.~T. Chen, J.~Ortiz,
  D.~DeTone, A.~Wang, S.~Anderson, {\em et~al.}, ``Theseus: A library for
  differentiable nonlinear optimization,'' {\em Advances in Neural Information
  Processing Systems}, vol.~35, pp.~3801--3818, 2022.

\end{thebibliography}
\vspace{0.1cm}
\begin{flushright}
  \scriptsize \framebox{\parbox{2.5in}{Government License: The
      submitted manuscript has been created by UChicago Argonne,
      LLC, Operator of Argonne National Laboratory (``Argonne").
      Argonne, a U.S. Department of Energy Office of Science
      laboratory, is operated under Contract
      No. DE-AC02-06CH11357.  The U.S. Government retains for
      itself, and others acting on its behalf, a paid-up
      nonexclusive, irrevocable worldwide license in said
      article to reproduce, prepare derivative works, distribute
      copies to the public, and perform publicly and display
      publicly, by or on behalf of the Government. The Department of Energy will provide public access to these results of federally sponsored research in accordance with the DOE Public Access Plan. http://energy.gov/downloads/doe-public-access-plan. }}
  \normalsize
\end{flushright}

\end{document}